\documentclass[letterpaper, 10 pt, conference]{ieeeconf}

\IEEEoverridecommandlockouts                              

\usepackage{algorithmic}
\usepackage{algorithm}
\usepackage{epstopdf}
\usepackage{amsmath}
\usepackage{amsfonts}
\usepackage{float}
\usepackage{breqn}
\usepackage{pgfplots}
\usepackage{subcaption}

\title{\LARGE \bf
Stochastic Active Discretizations for Accelerating Temporal Uncertainty Management of Gas Pipeline Loads
}

\author{Jake J. Harmon$^\star$, Svetlana Tokareva$^\star$, and Anatoly Zlotnik$^\star$
\thanks{This work was supported by the Laboratory Directed Research and Development (LDRD) program at
Los Alamos National Laboratory.  Research conducted at Los Alamos National Laboratory is done under the auspices of the National Nuclear Security Administration of the U.S. Department of Energy under Contract No. 89233218CNA000001. LA-UR 24-22636.}
\thanks{$^\star$\{harmon,tokareva,azlotnik\}@lanl.gov, \,\, Applied Mathematics and Plasma Physics Group, Theoretical Division, Los Alamos National Laboratory, PO Box 1663, Los Alamos, NM 87545, USA}}

\newtheorem{df}{Definition}

\begin{document}

\maketitle
\thispagestyle{empty}
\pagestyle{empty}

\begin{abstract}
We propose a predictor-corrector adaptive method for the simulation of hyperbolic partial differential equations (PDEs) on networks under general uncertainty in parameters, initial conditions, or boundary conditions.  The approach is based on the stochastic finite volume (SFV) framework that circumvents sampling schemes or simulation ensembles while also preserving fundamental properties, in particular hyperbolicity of the resulting systems and conservation of the discrete solutions.  The initial boundary value problem (IBVP) on a set of network-connected one-dimensional domains that represent a pipeline is represented using active discretization of the physical and stochastic spaces, and we evaluate the propagation of uncertainty through network nodes by solving a junction Riemann problem.  The adaptivity of our method in refining discretization based on error metrics enables computationally tractable evaluation of intertemporal uncertainty in order to support decisions about timing and quantity of pipeline operations to maximize delivery under transient and uncertain conditions.  We illustrate our computational method using simulations for a representative network. 
\end{abstract}

\section{Introduction} \label{secintro}
Discretized hyperbolic systems of conservation laws enable high-impact predictive computing across application domains, providing essential data in design and assessment workflows for critical physical systems and processes. Examples include Euler's equations and the so-called Shallow Water Equations. Many other problems in multiphysics, arising in part from conservative physical laws such as conservation of density or momentum exhibit hyperbolic character, such as the magnetohydrodynamic \cite{beg2013} or radiationhydrodynamic equations \cite{dorfi1998}. Such partial differential equation (PDE) systems exhibit challenging behavior and properties that significantly impact numerical analysis and discretization schemes, which complicates their use in practical simulation. In particular, solutions to these PDEs can form discontinuities in finite time, even from globally smooth initial conditions \cite{toro2009}.

A complicating factor in these physical models, many uncertainties, such as those originating from measurement or model errors, can inhibit predictive computing and predictive control, necessitating more sophisticated numerical schemes that furnish statistical quantities as opposed to only point estimates. Moreover, though in some sense secondary to ``physical'' uncertainties, discretization error can drive significant unreliability in estimating statistical quantities \cite{estep2009, butler2014}.  This need to control discretization error motivates self-guided, self-correcting methods that can efficiently extract uncertainty quantification (UQ) data with high confidence. Adaptive schemes can lead not only to improved accuracy, but also enhanced convergence rates and significant improvements to computational resource allocations \cite{berger1984, davis2016, deiterding2016, berger2021, harmon2022, harmon2022_sie, corrado2022, harmon2024adaptive}.

The rapid predictive simulation of uncertainty propagation in gas transport networks is a particularly compelling problem, because practical solutions to the associated hyperbolic PDE systems would enable event-triggered recourse actions for safe control.  For example, the ability to predict the probability of constraint violation hours in advance could indicate the likely need for corrective action before an event that would require remediation or curtailment of delivery to customers.   Uncertainties arise due to many sources in the setting of pipeline transport: gas supply status, compressor station conditions, gas withdrawals (customer demand), and prevailing weather and temperature, depending on the sophistication of the underlying model. Auxiliary sources of uncertainty, i.e., those independent of the physics such as economic cost, also significantly impact operations of gas pipeline networks. Incorporation of these uncertainties in optimal control workflows, however, presents its own challenges, particularly because discontinuities can propagate in the physical and stochastic spaces \cite{tokareva2014, mishra2016}.

Computing statistics via Monte Carlo integration, including in its multilevel forms, remains the dominant method for extending deterministic workflows to include uncertainty quantification \cite{giles2008, hoel2012, giles2015, elfverson2016, eigel2016, krumscheid2018, vanbarel2019}. While possessing favourable characteristics with respect to the dimension of the stochastic space and ease of implementation, Monte Carlo based methods converge slowly, requiring many thousands of samples. When the underlying physical model demands significant computational expenditure for even one realization, longer time-to-solution (even when supported by high-performance computing clusters) must be tolerated, or alternative methods must be sought. Such alternatives include generalized polynomial chaos (gPC) representations computed through stochastic Galerkin projection \cite{xiu2003, schlachter2020}. Unfortunately, gPC-based methods depend significantly on stochastic regularity. When applied to problems of hyperbolic character, where solutions generally lack smoothness, a gPC ansatz often suffers from spurious oscillations. Furthermore, the transformed stochastic Galerkin systems can, without special care, lose hyperbolicity, leading to questionable physical validity and predictive power.

Recently, a method based on hyperbolicity preserving, conservative, finite volume schemes demonstrated significant potential for enhancing the modeling of hyperbolic flows on networks with temporal uncertainty \cite{tokareva2024stochastic}. Furthermore, based on a similar formulation for conventional (i.e., non-networked) domains, an anisotropic method with enhanced convergence characteristics for UQ was introduced in \cite{harmon2024adaptive}. To facilitate UQ-informed predictive control of large-scale gas networks, we consider in this paper the extension of this adaptivity to hyperbolic flows on networks.

The remainder of this paper is organized as follows. In Section \ref{sec:main}, we outline the modeling of hyperbolic flows on networks. In Section \ref{sec:adaptivity}, we study a predictor-corrector scheme for automatically generating tuned discretizations that adapt to the modeling difficulty in the physical and stochastic spaces over the network. Finally, we consider an example network problem in Section \ref{sec:example}.

\section{Modeling Hyperbolic Flows on Networks}\label{sec:main}

Consider a network $\mathcal{N}$ of edges $\mathcal{E}$ and junctions $\mathcal{J}$. On each edge $i\in\mathcal{E},$ we assume the gas flows are governed by the PDE of form
\begin{equation}\label{eq:deterministic_cons_law}
\mathbf{u}_t + \nabla\cdot \textbf{F}(\mathbf{u}) = \mathbf{S}(\mathbf{u})\; \mathrm{in}\; \Omega_{i}\times\mathbb{I},\, \forall i\in\mathcal{E},   
\end{equation}
where $\Omega_i$ denotes the physical domain of the edge $i$, $\mathbb{I}$ denotes the temporal domain $[0,\, T]$ for $T\in \mathcal{R}_{>0}$, and $\mathbf{u}$ represents the state vector of conserved quantities. We assume that $\mathbf{F}$ is a hyperbolic flux, with its Jacobian having only real eigenvalues. The term $\mathbf{S}$ represents a source. In our notation the subscript $(\cdot )_t$ designates a partial derivative with respect to time $t$, while $\nabla$ applies to the physical coordinates. For clarity of notation, unless otherwise specified, we suppress references to the domains of the problem.

Many important problems have the form of \eqref{eq:deterministic_cons_law}, such as Euler's equations of gas dynamics with the flux
\begin{align}\label{eq:euler_flux}
    \mathbf{F}_{\mathrm{Euler}} &= \begin{bmatrix}
                                        \rho v \\
                                        \rho v^2 + p \\
                                        v(\rho E + p)
                                        \end{bmatrix},
\end{align} associated with the state vector
\begin{align}
    \mathbf{u} &= \begin{bmatrix}
                                        \rho \\
                                        \rho v \\
                                       \rho E
                                        \end{bmatrix},
\end{align}
where $\rho$, $v$, $p,\, E$ denote, respectively, the density, velocity, pressure, and total energy of the flow. Euler's equations also depends on an auxiliary equation of state (EOS), such as the ideal gas EOS.

For the modeling of flows on gas networks, we employ the following one-dimensional model \cite{misra2020, tokareva2024stochastic}
\begin{align}
    \rho_t + q_x &= 0, \label{eq:gas_1}\\
    q_t + a^2 \rho_x &= - \frac{\lambda}{2 D}\frac{q|q|}{\rho}, \label{eq:gas_2}
\end{align}
with $\rho$, as before, denoting the density and $q$ the mass flux. This model incorporates parameters of the network in the source term, namely the Darcy-Weisbach friction factor $\lambda$ and the pipe diameter $D$. The sound speed $a > 0$ is assumed constant. Of course, this system is supplemented by initial and boundary conditions such that it is well defined. Translating \eqref{eq:gas_1} and \eqref{eq:gas_2} to the form of \eqref{eq:deterministic_cons_law} leads to
\begin{align}
    \mathbf{u} = \begin{bmatrix}
                                        \rho \\
                                        q
                                        \end{bmatrix},
                                        \quad
                                        \mathbf{F} = \begin{bmatrix}
                                        q \\
                                        a^2 \rho
                                        \end{bmatrix},\quad
                                        \mathbf{S} = \begin{bmatrix}
                                        0 \\
                                        - \frac{\lambda}{2 D}\frac{q|q|}{\rho}
                                        \end{bmatrix}.
\end{align}

At the junctions $j\in\mathcal{J}$, compatibility conditions must be imposed due to the incident edges of the network. These conditions ensure, for example, conservation of mass as the gas flows through the network, eventually reaching customer withdrawals. More details on these conditions and the resulting junction Riemann problems may be found in our previous study \cite{tokareva2024stochastic}.

Until now, we have assumed the problem remains deterministic. Instead, consider uncertainty introduced in various parameters of the problem. Take $\left(\Omega_{\mathrm{stoch}}, \, \mathcal{F},\, \mathbb{P} \right)$, a probability space composed of a set $\Omega_{\mathrm{stoch}}$, a $\sigma$-algebra $\mathcal{F}$, and a probability measure $\mathbb{P}$ on $\mathcal{F}$.
Let $\mathbf{y} = \mathbf{y}(\omega)\in D_{\mathrm{stoch}}\subset\mathbb{R}^\ell$ measurable $\mathcal{F}$ be a random variable on this probability space and assume there exists for $\mathbf{y}$ a probability density $\mu$ such that $$\mathbb{P}(A) = \int_A \mu(\mathbf{y})\, d\mathbf{y}$$ for all $A\in \mathcal{F}$.

With this uncertainty, we obtain a new problem parameterized by the random variable $\mathbf{y}$ such that
\begin{equation}\label{eq:stochastic_cons_law_param}
    \mathbf{u}_t + \nabla\cdot\mathbf{F}(\mathbf{u}, \, \mathbf{y}) = 0, \; \mathbf{x}\in\Omega_{\mathrm{phys}},\, \mathbf{y}\in D_{\mathrm{stoch}},
\end{equation}
subject to
\begin{align}
 \mathbf{u}(\mathbf{x},0,\mathbf{y}) = \mathbf{u}_0(\mathbf{x},\mathbf{y}),&\quad \mathbf{x}\in \Omega_{\mathrm{phys}},\; \mathbf{y}\in D_{\mathrm{stoch}}, \label{eq:stochastic_ini_con_param} \\
    \mathbf{u}(\mathbf{x},t,\mathbf{y}) = \mathbf{u}_B(t,\mathbf{y}),&\quad \mathbf{x}\in\partial\Omega_{\mathrm{phys}},\; \mathbf{y}\in D_{\mathrm{stoch}}. \label{eq:stochastic_boundary_con_param}
\end{align}

In this setting, the state vector $\mathbf{u}$ is now a random variable, with its first moment given by
\begin{equation}
    \mathbb{E}[u] = \int_{\mathbb{R}^q} u\left(\cdot; \, \mathbf{y}\right)\mu\left(\mathbf{y}\right) d\mathbf{y} < \infty,
\end{equation} and higher moments computed similarly. Note that the probability density of $\mathbf{y}$ need not have finite support, and its specification depends on the \textit{principle of maximum entropy}; in particular, knowledge of minimum and maximum values alone implies a uniform distribution on $\mathbf{y}$, while knowledge of mean and variance implies Gaussian. 

We proceed according to our previous study on general active discretization methods for UQ in hyperbolic PDE systems \cite{harmon2024adaptive}, and discretize the stochastic space by a collection of cells $T_y$ such that
$$ \bar{D}_{\mathrm{stoch}} = \bigcup_{y} T_y, $$
i.e., the cells form a complete covering of $D_{\mathrm{stoch}}.$ Over each cell, or \textit{stochastic control volume}, the conditional expectation of $\mathbf{u}$ given $\mathbf{y}\in T_y$ is given by
\begin{equation}
\mathbb{E}[\mathbf{u}\,|\,\mathbf{y}\in T_y] = \frac{1}{P\left(\mathbf{y}\in T_y\right)}\int_{T_y} \mathbf{u}\mu(\mathbf{y})\,d\mathbf{y}. 
\end{equation} 
Note that the stochastic space is in some sense \textit{artificial}; its partitioning may be chosen freely, and therefore a structured mesh is often preferred.

Take a similar discretization of the physical space on each edge $i \in \mathcal{E}$ by cells $T_x$. Depending on the model, for instance $\Omega_{i}\subset\mathbb{R}^2$ or $\Omega_{i}\subset{\mathbb{R}}^3$, unstructured meshes may be preferred for this space. However, for the flows considered in this study, we assume physically one-dimensional pipes given the evidence in the literature for the predictive capability based on this representation \cite{gyrya2019, misra2020, tokareva2024stochastic}.

Rather than seek solutions pointwise, instead consider equality over the stochastic control volumes:
\begin{dmath}\label{eq:sfv_fv_scheme}
    \frac{1}{P\left(\mathbf{y}\in T_y\right)}\left( \int_{T_x} \int_{T_y} \mathbf{u}_t\mu(\mathbf{y})\,d\mathbf{x}d\mathbf{y} + \int_{T_x} \int_{T_y} \nabla\cdot\mathbf{F}(\mathbf{u}, \mathbf{y})\mu(\mathbf{y})\, d\mathbf{x}d\mathbf{y} \right) = \bar{\mathbf{S}}_y(\mathbf{u}),
\end{dmath}
where $\bar{\mathbf{S}}_y(\mathbf{u})$ denotes the conditional expectation of $\mathbf{S}$ over the stochastic control volume $T_y$. To simplify this expression for each edge, we introduce the measure $h_T$ of the combined computational cell $T=T_x\times T_y$ where

\begin{equation}
            h_T = |T_x||T_y| = \int_{T_x}\int_{T_y} \mu(\mathbf{y})\,d\mathbf{x}d\mathbf{y},
        \end{equation}
        and $|T_x|$ and $|T_y|$ denote the physical volume and the stochastic volume. We refer to the collection of such cells $T$ on edge $i$ by $\mathcal{T}_i$.

Supposing a piecewise constant ansatz on the domain of each edge $i\in \mathcal{E}$, let
        \begin{equation}
        {\mathbf{U}}_T = \frac{1}{h_T}\int_{T_x}\int_{T_y} \mathbf{u}(\mathbf{x},\,t,\,\mathbf{y})\mu(\mathbf{y})\,d\mathbf{x}d\mathbf{y}
        \end{equation}
        denote the average of the state vector over the computational cell $T$, i.e., an averaging over the physical \textit{and} stochastic spaces. This representation leads to an \textit{exact} PDE, now governing the averages:
         \begin{dmath}\label{eq:sfv_ODE}
            \frac{d{\mathbf{U}}_T}{dt} + \frac{1}{h_T}\int_{\partial T_x} \int_{T_y} \left( \mathbf{F}(\mathbf{u}, \mathbf{y})\cdot\hat{\mathbf{n}}\right) \mu(\mathbf{y})\, d\mathbf{x}d\mathbf{y} = \bar{\mathbf{S}},
        \end{dmath}
        where $\bar{\mathbf{S}}$ denotes the average of $\mathbf{S}$ over $T$.
        As in the conventional case, approximation error is introduced through the flux integral, due to the discontinuity of flows and, in this case, the need to estimate point values (specifically, at the boundaries of cells) from piecewise constant data.

        Regardless, for each computational cell $T$ associated with a portion of the physical edge and a control volume in the stochastic space, we can evolve $\mathbf{U}_T$ forward in time. Statistics of the state vector (or \textit{functionals} of the state vector), for example expectations or variances, may be readily obtained through trivial post-processing \cite{tokareva2024stochastic, harmon2024adaptive}.

\section{Adaptivity}\label{sec:adaptivity}

In Section \ref{sec:main} we introduced the underlying mechanism for obtaining approximations of gas flows under uncertainty. However, the central problem---\textit{choosing} how to partition for each edge $i\in\mathcal{E}$ the spaces $\Omega_i$ and $\Omega_{\mathrm{stoch}}$---remains unaddressed. Here, we develop the extension of adaptivity and error control to the case of hyperbolic flows on networks.

Now, both the networked and conventional cases can exhibit significant sensitivity to the treatment at the computational boundary. In the conventional case, typical conditions might include freeflow, slip, or no-slip conditions. In the \textit{networked} case, we have in addition to boundary conditions (for example due to injections or withdrawals), the compatibility conditions mentioned in Section \ref{sec:main} at every junction in the network. The treatment of these compatibility conditions is essential to the predictive quality of the modeled flows, and therefore any subsequent predictive control reliant upon the resulting uncertainty quantification provided.

Likewise, the treatment of adaptivity at these junctions is essential; ineffective constraints at the interfaces between each $\Omega_i$ and their representations of $D_{\mathrm{stoch}}$ can lead to significant performance bottlenecks. The generalization, therefore, of the recently introduced predictor-corrector scheme \cite{harmon2024adaptive} benefits from the following two-stage procedure:
\begin{enumerate}
    \item Predict discretization error for every cell $T$ on every edge $i\in\mathcal{E}$ \textit{using local data}
    \item For each $i,\, j \in \mathcal{E}\left(\mathcal{J}_k\right)$, where $\mathcal{E}\left(\mathcal{J}_k\right)$ denotes the set of edges associated with the $k$th junction, propagate refinement directives via discretization compatibility requirements
\end{enumerate}

We outline Stage 1 below, and summarize the details in Algorithm \ref{alg:adaptivity}. The constraint of \textit{local data} is of some consequence, particularly for distributed memory environments, where communications between edges that share a junction but belong to different nodes on a cluster should be minimized. 
To proceed, we invoke the notion of a \textit{reconstructor}.
\begin{df}[Reconstruction, Def. (3.7, \cite{harmon2022})] A reconstructor $\mathbb{B}_{U(\mathcal{T}_i)}(T; \xi):\mathbb{V}(\mathcal{T}_i)\rightarrow\mathbb{R}^p$ maps a solution $\mathbf{U}(\mathcal{T}_i)\in\mathbb{V}(\mathcal{T}_i)$ associated with degrees of freedom on $\mathcal{T}$ to point values $\xi\in T$.
\end{df}

Note that this notion of reconstruction is local in the sense that \textit{only data on an edge is used to reconstruct point values on that same edge}. Furthermore, though the underlying approximation, namely the degrees of freedom, may be piecewise constant, the reconstruction from this data may be higher-order, e.g., piecewise polynomial. Crucially, from the same data, we may obtain multiple reconstructions. In our recent study \cite{harmon2024adaptive}, this property is exploited to estimate discretization error by predicting error at future time step from the error between an enriched reconstructor $\mathbb{B}^{\mathrm{H}}$ and a reduced reconstructor $\mathbb{B}^{\mathrm{L}}$. In other words, we obtain a predictor $\mathcal{P}$ such that
\begin{equation}
    \eta_T \leftarrow \mathcal{P}\left(\mathbf{U}(\mathcal{T}_i),\, T, \,\mathbb{B}^{\mathrm{H}},\, \mathbb{B}^{\mathrm{L}} \right)
\end{equation}
delivers a refinement indicator $\eta_T$ for the cell $T$. This predictor is applied for all edges, resulting in highlighting of sufficient, insufficient, and inefficient regions of the computational domain. The predictor $\mathcal{P}$ is computed until satisfaction of desired accuracy tolerance, leading to a so-called predictor-corrector; until the predicted and observed results match within a desire tolerance, refinement must continue.
\begin{figure}[!htb]
        \centering
       \begin{subfigure}[b]{0.75\linewidth}{
       \includegraphics[]{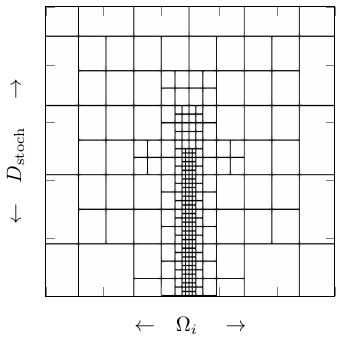}
       }
        \caption{Isotropic}
        \label{fig:sub:iso}
\end{subfigure}\\
\begin{subfigure}[b]{0.75\linewidth}{
       \includegraphics[]{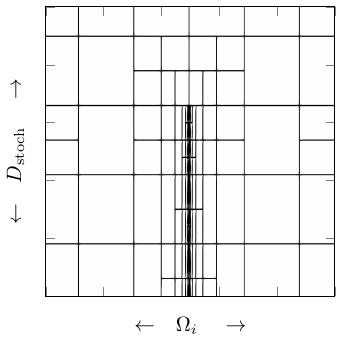}
}
        \caption{Anisotropic}
        \label{fig:sub:aniso}
\end{subfigure}
\caption{Refinement paradigms on a given edge $i\in \mathcal{E}$ of the gas network under uncertainty. (a) Isotropic refinements in the physical and stochastic spaces. (b) Anisotropic refinements in the physical and stochastic spaces. Anisotropic refinements can capture non-smooth flows with targeted, directional allocations of computational resources for more efficient UQ-informed predictive control.}
        \label{fig:iso_vs_aniso}
\end{figure}

Yet, as in the conventional case, knowing \textit{where} to refine, does not imply \textit{how} to refine. For this, we recall the smoothness based flag introduced in \cite{harmon2024adaptive} for deciding in which directions to partition the combined physical and stochastic spaces, 
\begin{dmath}\label{eq:smoothness_indicator}
    \beta_{\mathcal{O}} = \sum_{i = 1}^{\ell} \int_{T_{\mathcal{O}}} (\Delta \xi)^{2i - 1} \left(  \frac{\partial^i \mathbb{B}^{1\mathrm{D}}_{U(\mathcal{T})}(T;\, \xi,\, \mathcal{O})}{\partial^i \xi}  \right)^2\, d\xi,
\end{dmath}
where $\mathcal{O}$ represents a target direction, $T_{\mathcal{O}}$ a slice of the cell $T$ in the $\mathcal{O}$ direction, $\mathbb{B}^{1\mathrm{D}}_{U(\mathcal{T})}$ a one-dimensional reconstruction in the same direction, and $\ell$ denotes a smoothness degree. Typically, we let $\ell = 1$ to decide the importance of the directions in the computational domain.

In contrast to \textit{isotropic} refinements, which ignore directionality of the flows and their modeling difficulty in the physical and stochastic directions, \textit{anisotropic} refinements distinguish between the spaces and their interactions. We illustrate this distinction for a single one-dimensional edge $i\in\mathcal{E}$ with an associated one-dimensional stochastic space in Fig. \ref{fig:iso_vs_aniso}.

\begin{figure}[!tb]
\centering
        \hspace{0.5cm}
        \includegraphics[]{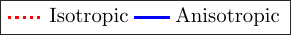}
{
\begin{center}
    \hspace{-0.75cm}
\includegraphics[]{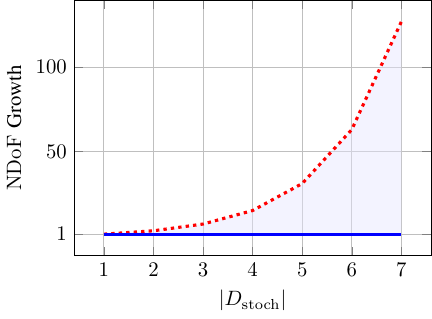}
\end{center}
}
\caption{Minimum growth in the number of degrees of freedom (NDoFs) with respect to the dimension of the stochastic space upon refinement of a stochastic control volume. Anisotropic refinements have flexibility in allocating new resources, whereas isotropic refinements are constrained to exponential growth.}
        \label{fig:dof_growth}
\end{figure}

\textit{Qualitatively}, anisotropic refinements deliver full refinement flexibility. \textit{Quantitatively}, refinements of anisotropic character enable significant enhancements in the growth of the number of degrees of freedom (NDoFs), as illustrated in Fig. \ref{fig:dof_growth}.

Finally, with a D\"{o}rfler-like selection criterion on $\beta_{\mathcal{O}},$
\begin{equation}\label{eq:anisotropy_decision} 
\beta_{\mathcal{O}} > \varepsilon_{\mathrm{aniso}} \sum^{d+q}_{i=1} \beta_{i},\quad \varepsilon_{\mathrm{aniso}} < \frac{\max_i{\beta_i}}{\sum^{d+q}_{i=1} \beta_{i}}, 
\end{equation} 
we obtain a list of refinement instructions $\left\{ \mathcal{R} \right\}_{i\in\mathcal{E}}$ for each edge in the network.

\begin{algorithm}
\caption{Unsteady Enriched-Reduced $\mathbb{B}$ Adaptivity \cite{harmon2024adaptive}}
\label{alg:adaptivity}
\begin{algorithmic}[1]
\STATE{Define \{$\varepsilon, \, \varepsilon_{\mathrm{aniso}},\, t_n,\, \mathbb{B}^{\mathrm{H}},\, \mathbb{B}^{\mathrm{L}},\,\theta,\, \mathcal{T}$\}}
\FORALL{$i\in\mathcal{E}$}
\STATE{${\eta} := \infty$}
\WHILE{${\eta} > \varepsilon$}
\STATE{${\eta} = 0$}
\FORALL{$T\in\mathcal{T}$}
\STATE\label{line3}{Perform reconstruction and flux integration for an enriched-reduced reconstruction pair $\mathbb{B}^{\mathrm{H}}$ / $\mathbb{B}^{\mathrm{L}}$}
\STATE{Estimate a local time update error $\eta_T$ according to Prop. (4.5, \cite{harmon2024adaptive})}
\IF{$\eta_T > \varepsilon$}
\STATE{Compute smoothness indicators via \eqref{eq:smoothness_indicator}}
\STATE{Assign refinement directions via \eqref{eq:anisotropy_decision}}
\STATE{Mark for refinement}
\ELSIF{$\eta_T < \theta \varepsilon$}
\STATE{Mark for coarsening}
\ENDIF
\STATE{${\eta} = \max{\left( \eta,\, \eta_T \right)}$}
\ENDFOR
\ENDWHILE
\ENDFOR
\RETURN $\left\{ \mathcal{R} \right\}_{i\in\mathcal{E}}$
\end{algorithmic}
\end{algorithm}
Without the additional compatibility conditions arising from the network architecture of the problems under consideration, this process, summarized in Algorithm \ref{alg:adaptivity}, is sufficient to automatically control the error in the evolution of hyperbolic flows. The extension, however, of this concept to networks requires discretization compatibility requirements as part of a second refinement stage.

Note that each edge $i\in\mathcal{E}$, in addition to its ``natural'' ownership of the discretization of $\Omega_i$, possesses its own \textit{realization} of $D_{\mathrm{stoch}}$. As the character of flows need not be uniform, particularly for large networks, shared discretization or representation of the stochastic space across the problem domain augments computational expenditure due to inefficient allocation of degrees of freedom \textit{and} significant communication bottlenecks needed to reconcile refinement requests across the entire controlled network. 

Instead, we augment a subset of the instructions $\left\{ \mathcal{R} \right\}_{i\in\mathcal{E}}$ to enforce interfacing matching conditions. The concept of \textit{flux 1-irregularity} introduced in \cite{harmon2024adaptive} is insufficiently strict for resolving the flux integrals at the junctions. Hence, for each $i,\, j \in \mathcal{E}\left(\mathcal{J}_k\right)$, letting $\left\{T_{\partial \Omega_i}\right\}$ and $\left\{T_{\partial \Omega_j}\right\}$ denote the sets of cells incident at the junction, we stipulate that if there exists a refinement instruction in $\mathcal{R}_i$ associated with $\left\{T_{\partial \Omega_i}\right\}$, an equivalent instruction in $\mathcal{R}_j$ must exist. In other words, the representation of the stochastic space at the junction is constrained to uniformity, with free variation elsewhere in the computational domain.

This process is repeated, with resources allocated and de-allocated as necessary, until all quantities and statistics required for predictive control workflows are obtained.

\section{Case Study}\label{sec:example}

As an illustrative numerical study, we consider a previously defined and examined IBVP for a small gas pipeline test network with intertemporal uncertainty \cite{tokareva2024stochastic}, which consists of five junctions and five edges. The full specification is detailed in the previous study (Section 4.2, \cite{tokareva2024stochastic}), and for the sake of brevity we forgo its replication here. We do however note that uncertainty is imposed at the terminal point of the network, with random variation in the customer withdrawal rate at the fifth node. The uncertainty is temporal in nature, corresponding to stochasticity in the initial time of an augmentation in gas consumption. Explicitly, we take the withdrawal on this terminal node as
\begin{equation}\label{eq:withdrawal}
    d(t) = \begin{cases}
      d_1(t), & \text{if }t < \tau_p\\
      \frac{d_2 - d_1}{\tau_p^1 - \tau_p}(t-\tau_p) + d_1(t), & \text{if }t\in(\tau_p,\,\tau_p^1)\\
      d_2(t), & \text{if }t\in(\tau_p^1,\, \tau_p^2)\\
      \frac{d_1-d_2}{\tau_p^3 - \tau_p^2}(t-\tau_p^2) + d_2(t), & \text{if }t\in(\tau_p^2,\,\tau_p^3)\\
      d_1(t), & \text{if }t> \tau_p^3
    \end{cases}    
\end{equation}
where $\tau_p$ denotes an uncertain time of withdrawal, and $(\tau_p - \tau_p^1),\, (\tau_p - \tau_p^2),$ and $(\tau_p - \tau_p^3)$ are deterministic event intervals. The withdrawal functions $d_1$ and $d_2$ are deterministic to isolate the uncertainty propagation that results from the random time of the event alone. See \cite{tokareva2024stochastic} for their specification.

For the temporal domain, we consider a 24 hour period. At initialization, according to the specification in \cite{tokareva2024stochastic}, we assume the flows are deterministic. Uncertainty propagates throughout the network as a result of an uncertain time of withdrawal according to \eqref{eq:withdrawal} during this 24 hour period. We further suppose $\tau_p\sim \mathcal{U}(4, 12),$ but any other distribution on $\tau_p$ could be chosen.

Along with the networked unsteady enriched-reduced $\mathbb{B}$ adaptivity we proposed in Section \ref{sec:adaptivity}, we employ third-order strong stability preserving (SSP) time-integration for evolving the discretized system of equations.

\begin{figure}[!tb]
    \centering
    \hspace{-0.5cm}
    \includegraphics[]{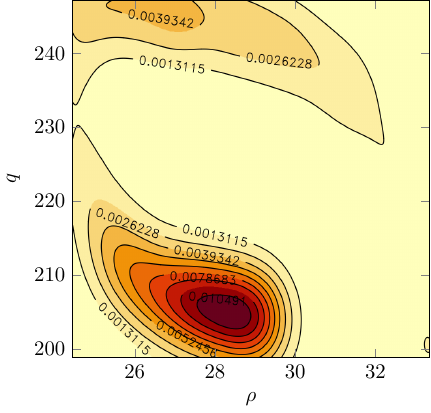}
    \caption{The joint push-forward probability density function of the density $\rho$ and mass flux $q$ in the middle of terminal edge of the test network at hour 12.}
    \label{fig:push_forward_density}
\end{figure}

Accurate statistics of the flow quantities are essential for predictive control. We therefore consider the push-forward probability densities of $\rho$ and $q$. In Fig. \ref{fig:push_forward_density}, we consider the joint push-forward probability density function in the middle of the terminal edge of the network at hour 12. By this time in the evolution of the network, uncertainty in the withdrawal time has propagated throughout the networked, leading to the rendered statistical relationship between the components of the state vector.

\begin{figure}[!htb]
    \begin{centering}
        \hspace{1.0cm}
        \includegraphics[]{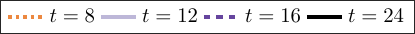}
     \end{centering}
       \begin{center}
        \hspace{-0.75cm}
       \includegraphics[]{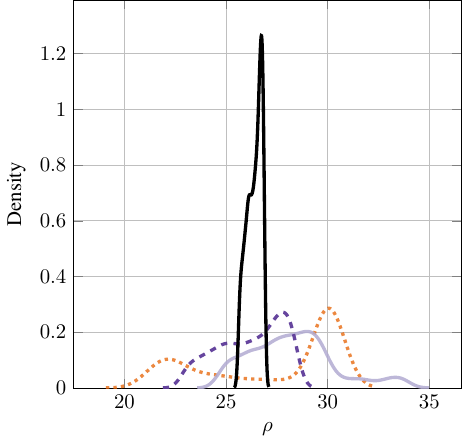}
        \end{center}
\caption{Push-forward probability density function for the density $\rho$ in the middle of the terminal edge at selected simulation times $t$ in hours. Prior to hour 4, the flows are deterministic. As the maximum possible timing of the uncertain withdrawal (12 hours) is passed, the network stabilizes and the uncertainty narrows from complicated, multi-modal distributions.}
        \label{fig:uq_prop_density}
\end{figure}

Next, we consider in Fig. \ref{fig:uq_prop_density} the marginal push-forward densities at a series of times to illustrate the propagation of uncertainty in the network and the level of UQ afforded by our method. At initialization until hour 4, the flows on the network are deterministic. Due to the distribution on $\tau_p$, by hour 8, the push-forward probability density on $\rho$ becomes multi-modal. After hour 12 (the maximum possible time of the withdrawal event), the uncertainty again narrows, leading to the distribution on $\rho$ seen for the terminal time. Our method efficiently exposes this highly varied, temporal fluctuation in uncertainty for enhanced predictive control.

\section{Conclusion}\label{sec:conclusion}
We constructed the extension of predictor-corrector adaptivity to the case of hyperbolic flows on networks. By treating the modeling difficulty throughout the network through active discretization of the stochastic and physical spaces, new, more efficient computations can be performed. As part of a predictive control workflow, enhanced model quality, both in terms of accuracy and efficiency, combined with uncertainty quantification can alleviate computational burdens associated with probabilistically constrained optimization objectives. Future work will involve the coupling of this enhanced gas network simulator with stochastic optimization for control of realistic networks.
 \bibliographystyle{unsrt}
 \bibliography{main}

\end{document}